\documentclass[10pt]{article}

\usepackage{amssymb,amsmath}
\usepackage{nth}

\usepackage{url}

\newtheorem{theorem}{Theorem}[section]
\newtheorem{corollary}[theorem]{Corollary}

\newtheorem{lemma}[theorem]{Lemma}
\newtheorem{conjecture}[theorem]{Conjecture}
\newtheorem{example}[theorem]{Example}
\newtheorem{proposition}[theorem]{Proposition}

\numberwithin{equation}{section}

\newcommand{\Qed}{\hfill $\Box$ \medskip}

\begin{document}

\title{A note on $m$-factorizations of complete multigraphs 
arising from designs} 
\author{Gy. Kiss$^{\dag}$
\thanks{The research was supported by the Mexican-Hungarian Intergovernmental
Scientific and Technological Cooperation Project, 
Grant No. T\'ET 10-1-2011-0471.
\newline \indent
$^{\dag}$Author was supported by 
the Hungarian National Foundation for Scientific Research, Grant No. K 81310.}
\, and \, C. Rubio-Montiel$^*$}
\date{}

\maketitle


\begin{abstract}
Some new infinite families of simple, indecomposable $m$-factorizations 
of the complete multigraph $\lambda K_v$ are presented.
Most of the constructions come from finite geometries.
\end{abstract}

\section{Introduction}

The \emph{complete multigraph} $\lambda K_{v}$ has $v$ vertices 
and $\lambda$ edges
joining each pair of vertices. An \emph{$m$-factor} of the complete multigraph
$\lambda K_{v}$ is a set of pairwise vertex-disjoint $m$-regular subgraphs,
which induce a partition of the vertices. An \emph{$m$-factorization} 
of $\lambda K_{v}$ is a set of pairwise edge-disjoint
$m$-factors such that these $m$-factors induce a partition of the edges. An
$m$-factorization is called \emph{simple} if 
the $m$-factors are pairwise distinct.  
Furthermore, an $m$-factorization of $\lambda K_{v}$ is
\emph{decomposable} if there exist positive 
integers $\mu_{1}$ and $\mu_{2}$
such that $\mu_{1}+\mu_{2}=\lambda$ and the factorization is the union
of the $m$-factorizations $\mu_{1}K_{v}$ and $\mu_{2}K_{v}$, otherwise
it is called \emph{indecomposable}. There is no direct
correspondence between simplicity and indecomposability. 

Many papers deal with $m$-factorizations of graphs and multigraphs. 
This is an interesting problem in its own right, but it is motivated 
by several applications, too. In particular 
if $m=1,$ then a one-factorization of $K_{v}$ corresponds to a schedule of a
round robin tournament. For
a comprehensive survey on one-factorizations we refer to \cite{MR1433596}.
A special case of $2$-factorizations is the famous \emph{Oberwolfach problem},
see e.g. \cite{alsp, 2-1}.
Several authors investigated $3$-factorizations of $\lambda K_v$ with a 
certain 
automorphism group, see e.g. \cite{3-1, 3-2}. In general, decompositions of 
$\lambda K_v$ 
is also a widely studied problem, see e.g. \cite{m-2, m-3, m-4, m-1}. 
As $m$ increases, the structure of an 
arbitrary $m$-factor of $\lambda K_v$ can be much more complicated and 
the existence
problem becomes much more difficult. In this paper we restrict ourselves to 
construct factorizations in which all factors are regular graphs of 
degree $m$ whose connected components are complete graphs on $(m+1)$
vertices. In the case $m=1$ 
an indecomposable one-factorization of
$\lambda K_{2n}$ is denoted by $\mathrm{IOF}(2n,\lambda)$. 
Only a few conditions on the parameters are known: 
if $\mathrm{IOF}(2n,\lambda)$ exists, then 
$\lambda < 1\cdot3\cdot...\cdot(2n-3)$ \cite{MR1020644}; each
$\mathrm{IOF}(2n,\lambda)$ can be embedded in a 
simple $\mathrm{IOF}(2s,\lambda)$, 
provided that $\lambda < 2n < s$ \cite{MR802723}.
Six infinite classes of indecomposable one-factorizations have been
constructed so far, namely a simple $\mathrm{IOF}(2n,n-1)$ when $2n-1$ 
is a prime \cite{MR802723}, $\mathrm{IOF}(2(\lambda+p),\lambda)$ 
where $\lambda >2$ and $p$ is the
smallest prime wich does not divide $\lambda$ \cite{MR1140569} (an improvement
of this result can be found 
in \cite{MR1822868}), a simple $\mathrm{IOF}(2^h+2,2)$ where
$h$ is a positive integer \cite{MR1821979}, $\mathrm{IOF}(q^2+1,q-1)$ 
where $q$ is an odd prime number \cite{MR1805305}, a simple 
$\mathrm{IOF}(q^2+1,q+1)$ for any odd prime power $q$ \cite{MR1883515}, 
and a simple $\mathrm{IOF}(q^2,q)$ for any even prime power $q$
\cite{MR1883515}. Most of these constructions arise from finite geometry.

The aim of this paper is to construct 
new simple and indecomposable $m$-factorizations of $\lambda
K_{v}$ for different values of $m$, $\lambda$ 
and $v$. In Section \ref{Preliminaries} we recall the basic combinatorial
properties of designs and the geometric properties of finite affine and 
projective spaces. We also describe a general construction
method of $m$-factorizations which is based on spreads of block designs. 
In Sections \ref{Affine} and \ref{Proj} affine spaces and 
projective spaces, respectively, 
are the key objects. We present several new multigraph 
factorizations using subspaces, subgeometries and other configurations 
of these structures.

\section{Preliminaries}
\label{Preliminaries}
\renewcommand{\labelenumi}{(\alph{enumi})}

In this section we collect some concepts and results from design theory. 
For a detailed introduction to block designs we refer to \cite{cvl}.

\subsection{Designs}
Let $v$, $b$, $k$, $r$ and $\lambda$ be positive integers with $v>1$. 
Let $D=(\mathcal{P},\mathcal{B},\mathrm{I} )$ be a triple
consisting of a set
$\mathcal{P}$ of $v$ distinct objects, called points of $D$, a set
$\mathcal{B}$ of $b$ distinct objects, called blocks of $D$, and an incidence
relation $\mathrm{I}$, a subset of $\mathcal{P}\times \mathcal{B}$. 
We say that $x$ is incident with
$y$ (or $y$ is incident with $x$) if and only if the ordered pair $(x,y)$ is
in $\mathrm{I}$. $D$ is called a \emph{$2-(v,b,k,r,\lambda )$ design}
if it satisfies the following axioms.
\begin{enumerate}
\item 
Each block of $D$ is incident with exactly $k$ distinct points of $D$. 
\item 
Each point of $D$ is incident with exactly $r$ distinct blocks of $D$. 
\item 
If $x$ and $y$ are distinct points of $D,$ then there are exactly $\lambda$ 
blocks of $D$ incident with both $x$ and $y$. 
\end{enumerate}

A  \emph{$2-(v,b,k,r,\lambda )$} design is called a balanced incomplete
block design and is denoted by  \emph{$(v,k,\lambda )$}-design, too.
The parameters of a \emph{$2-(v,b,k,r,\lambda )$} design are not all
independent.  The two basic equations connecting them
are the following:
\begin{equation} \label{1.2}
vr=bk \quad \mathrm{and} \quad r(k-1)=\lambda(v-1).
\end{equation}
These necessary conditions are not sufficient, 
for example no $2-(43,43,7,7,1)$ design exists.

\subsection{Resolvability}

A \emph{resolution class} (or, a parallel class) of a 
\emph{$(v,k,\lambda )$}-design 
is a partition of the point-set of the design into blocks.
In general, an \emph{$f$-resolution class} 
of a design is a collection of blocks, which together contain
every point of the design exactly $f$ times.
A \emph{resolution} of a design is a partition of the block-set of the
design into $r$ resolution classes. A \emph{$(v,k,\lambda )$}-design 
with a resolution is called \emph{resolvable}.

Necessary conditions for the existence of a resolvable 
\emph{$(v,k,\lambda )$}-design are 
$\lambda (v-1)\equiv 0$ (mod $(k-1)$),
$v\equiv 0$ (mod $k$) and
$b\geq v+r-1,$ (see \cite{bose}).

Let $D=(\mathcal{P},\mathcal{B},\mathrm{I} )$ be a
\emph{$(v,k,\lambda )$}-design, where
$ \mathcal{P}=\{p_1,p_2,\dots, p_v\} $ is the set of its points 
and $\mathcal{B}=\{B_1,B_2,\dots, B_b\} $ is the set of its blocks. 
Identify the points of $D$
with the vertices of the complete multigraph $\lambda K_v$. Then in the
natural way, the set of points of each block of $D$ induces in $\lambda K_v$ a
subgraph isomorphic to $K_{k}$.
For $B_i\in \mathcal{B}$, let $G_i$ be the subgraph of $\lambda K_v$
induced by $B_i$. Then it follows from the properties of $D$ that a 
resolution class of $D$ gives a $(k-1)$-factor of $\lambda K_v$ and 
a resolution of $D$ gives a $(k-1)$-factorization of $\lambda K_v$. 
Hence we get the following well-known fact.

\begin{lemma}[Basic Construction]
\label{basic}
The existence a resolvable \emph{$(v,k,\lambda )$}-design is equivalent to
the existence of a $(k-1)$-factorization of the
complete multigraph $\lambda K_v$.   
\end{lemma}

\subsection{Projective and affine spaces}

Most of our factorizations come from finite geometries. In this subsection we
collect the 
basic properties of these objects. For a more detailed introduction we refer
to the book of Hirschfeld \cite{MR554919}.

Let $V_{n+1}$ be an $(n+1)$-dimensional vector space over the finite field of
$q$ elements, $\mathrm{GF}(q).$
The \emph{$n$-dimensional projective space} $\mathrm{PG}(n,q)$ is the geometry
whose 
$k$-dimensional subspaces for 
$k=0,1,\ldots ,n$ are the $(k+1)$-dimensional subspaces of 
$V_{n+1}$ with the zero deleted. 
A $k$-dimensional subspace of $\mathrm{PG}(n,q)$ is called $k$-space.
In particular subspaces of dimension zero, one and two are respectively 
a \emph{point,} a \emph{line} and a \emph{plane,} while a subspace of
dimension $n-1$ is called 
a \emph{hyperplane.}

The relation $\sim $
$${\mathbf x}\sim {\mathbf y}\Leftrightarrow \exists \,  
0\neq \alpha \in \mathrm{GF}(q) \, : \, 
{\mathbf x}=\alpha  {\mathbf y}$$
is an equivalence relation on the elements of $V_{n+1}\setminus \mathbf{0}$ 
whose equivalence classes are the points of $\mathrm{PG}(n,q).$
Let ${\mathbf v}=(v_0,v_1,\ldots ,v_n)$ be 
a vector in $V_{n+1}\setminus \mathbf{0}.$ 
The equivalence class of ${\mathbf v}$
is denoted by $[{\mathbf v}].$
The homogeneous coordinates of the point represented by $[{\mathbf v}]$ 
are $(v_0:v_1:\ldots :v_n).$ Hence two $(n+1)$-tuples $(x_0:x_1:\ldots :x_n)$
and $(y_0:y_1:\ldots :y_n)$ represent the same point of $\mathrm{PG}(n,q)$ if
and only if there exists $0\neq \alpha \in \mathrm{GF}(q)$ such that 
$x_i=\alpha y_i$ holds for $i=0,1,\ldots ,n.$ 

A $k$-space contains those points whose representing vectors ${\mathbf x}$ 
satisfy the equation ${\mathbf x}A={\mathbf 0},$ where $A$ is an 
$(n+1)\times (n-k)$ matrix of rank $n-k$ with entries in $\mathrm{GF}(q).$ 
In particular a hyperplane contains those points whose homogeneous 
coordinates $(x_0:x_1:\ldots :x_n)$ satisfy a linear equation  
$$u_0x_0+u_1x_1+\dots +u_nx_n=0$$
where $u_i\in \mathrm{GF}(q)$ and $(u_0,u_1,\ldots ,u_n)\neq {\mathbf 0}.$    

The basic combinatorial properties of $\mathrm{PG}(n,q)$ can be described 
by the $q$-nomial coefficients. $\genfrac[]{0pt}{2}{n}{k} _q$ equals to the
number of $k$-dimensional subspaces in an $n$-dimensional vector space over 
$\mathrm{GF}(q),$ hence it is defined as
$$ \genfrac[]{0pt}{0}{n}{k} _q := 
\frac{(q^n-1)(q^n-q)\ldots (q^n-q^{k-1})}{(q^k-1)(q^k-q)\ldots (q^k-q^{k-1})}.$$
The proof of the following proposition is straightforward. 
\begin{proposition}
\label{combprop}
$ $
\begin{itemize}
\item
The number of $k$-dimensional subspaces in $\mathrm{PG}(n,q)$ is 
$\genfrac[]{0pt}{1}{n+1}{k+1} _q.$
\item
The number of $k$-dimensional subspaces of $\mathrm{PG}(n,q)$ through a given 
$d$-dimensional ($d\leq k$) subspace in PG$(n,q)$ is 
$\genfrac[]{0pt}{1}{n-d}{k-d} _q.$ 
\item
In particular the number 
of $k$-dimensional subspaces of $\mathrm{PG}(n,q)$ through two distinct points 
in $\mathrm{PG}(n,q)$ is 
$\genfrac[]{0pt}{1}{n-1}{k-1} _q.$
\end{itemize}
\end{proposition}

If ${\mathcal H}_{\infty }$ is any hyperplane of $\mathrm{PG}(n,q),$ then
the \emph{$n$-dimensional affine space over $\mathrm{GF}(q)$} is 
$\mathrm{AG}(n,q)=\mathrm{PG}(n,q)\setminus {\mathcal H}_{\infty }.$
The subspaces of $\mathrm{AG}(n,q)$ are the subspaces of $\mathrm{PG}(n,q)$
with the points of ${\mathcal H}_{\infty }$ deleted in each case. 
The hyperplane ${\mathcal H}_{\infty }$ is called the 
\emph{hyperplane at infinity of $\mathrm{AG}(n,q),$} 
and for $k=0,1,\ldots, n-2$ the $k$-dimensional subspaces 
in ${\mathcal H}_{\infty }$ are called the 
\emph{$k$-spaces at infinity of $\mathrm{AG}(n,q).$} 
Let $1<d<n$ be an integer. Two $d$-spaces of $\mathrm{AG}(n,q)$ are
called \emph{parallel}, if the corresponding $d$-spaces of $\mathrm{PG}(n,q)$
intersect ${\mathcal H}_{\infty }$ in the same $(d-1)$-space. The parallelism
is an equivalence relation on the set of $d$-spaces of $\mathrm{AG}(n,q)$. As
a straightforward corollary of Proposition \ref{combprop} we get the following.

\begin{proposition}
\label{parallelnumber}
In $\mathrm{AG}(n,q)$ each equivalence class of parallel $d$-spaces 
contains $q^{n-d}$ subspaces. 
\end{proposition}

Projective and affine spaces provide examples of designs.

\begin{example}
\label{ex-proj}
Let $i<n$ be positive integers.
The projective space $\mathrm{PG}(n,q)$ 
can be considered as a 2-design $D=(\mathcal{P},\mathcal{B},\mathrm{I} ),$
where $\mathcal{P}$ is the set of points of $\mathrm{PG}(n,q)$, $\mathcal{B}$
is the set of $i$-spaces of $\mathrm{PG}(n,q)$ and $\mathrm{I} $ is the set 
theoretical inclusion. The parameters of $D$ are $v=\tfrac{q^{n+1}-1}{q-1}$, 
$b=\genfrac[]{0pt}{1}{n+1}{i+1}_q$,
$k=\tfrac{q^{i+1}-1}{q-1}$, $r=\genfrac[]{0pt}{1}{n}{i}_q$ and
$\lambda =\genfrac[]{0pt}{1}{n-1}{i-1}_q$. 
\end{example}

\begin{example}
\label{ex-affine}
Let $i<n$ be positive integers.
The affine space $\mathrm{AG}(n,q)$ 
can be considered as a 2-design $D=(\mathcal{P},\mathcal{B},\mathrm{I} ),$
where $\mathcal{P}$ is the set of points of $\mathrm{AG}(n,q)$, $\mathcal{B}$
is the set of $i$-spaces of $\mathrm{AG}(n,q)$ and $\mathrm{I} $ is the set 
theoretical inclusion. The parameters of $D$ are 
$v=q^{n}$, $b=q^{n-i}\genfrac[]{0pt}{1}{n}{i}_q$,
$k=q^i$, $r=\genfrac[]{0pt}{1}{n}{i}_q$ and
$\lambda =\genfrac[]{0pt}{1}{n-1}{i-1}_q$. 
\end{example}

\noindent
In the rest of this paper Examples \ref{ex-proj} 
and \ref{ex-affine} will be denoted
by $\mathrm{PG}^{(i)}(n,q)$ and by $\mathrm{AG}^{(i)}(n,q),$ respectively.
We will use the terminology from geometry.
An \emph{$i$-spread}, $\mathcal{S}^i$, of $\mathrm{PG}(n,q)$ 
(or of $\mathrm{AG}(n,q)$) is a set of pairwise disjoint $i$-dimensional 
subspaces which gives a partition of the 
points of the geometry. In general, an \emph{$f$-fold $i$-spread}, 
$\mathcal{S}^i_f ,$ is a set of $i$-dimensional subspaces 
such that every point of the geometry is contained in exactly $f$ subspaces 
of $\mathcal{S}^i_f.$ An \emph{$i$-packing}, $\mathcal{P}^i$, of 
$\mathrm{PG}(n,q)$ (or of $\mathrm{AG}(n,q)$) is a set of spreads such that 
each $i$-dimensional subspace of the geometry is contained in 
exactly one of the spreads in $\mathcal{P}^i$, i.e., the spreads give 
a partition of the $i$-dimensional subspaces of the geometry.
The $i$-spreads, $f$-fold $i$-spreads and $i$-packings induce a resolution 
class, an $f$-resolution class and a resolution in $\mathrm{PG}^{(i)}(n,q)$ 
(or in $\mathrm{AG}^{(i)}(n,q)$), respectively. 

It is easy to construct spreads and packings in $\mathrm{AG}^{(i)}(n,q),$
because each parallel class of $i$-spaces is an $i$-spread. The situation 
is much more complicated in $\mathrm{PG}^{(i)}(n,q).$ There are only a few
constructions of spreads. The following theorem summarizes the known  
existence conditions.

\begin{theorem}[\cite{MR554919}, Theorems 4.1 and 4.16]
\label{spreadexists}
$ $
\begin{itemize}
\item
There exists an $i$-spread in $\mathrm{PG}^{(i)}(n,q)$ if and only if 
$(i+1)|(n+1).$
\item
Suppose that $i,l$ and $n$ are positive integers such 
that $(l+1)|\gcd (i+1,n+1).$ 
Then there exists an $f$-fold $i$-spread in $\mathrm{PG}^{(i)}(n,q),$ 
where $f=(q^{i+1}-1)/(q^{l+1}-1).$
\end{itemize} 
\end{theorem} 

There exist several different $1$-spreads (line spreads) in $\mathrm{PG}^{(1)}(3,q).$ 
We briefly mention two types. Let $\ell _1, \ell _2$ and  $\ell _3$ be three skew lines in $\mathrm{PG}(3,q).$ 
The set of the $q+1$ transversals of $\ell _1, \ell _2$ and  $\ell _3$ is called \emph{regulus}
and it is denoted by ${\cal R}(\ell _1, \ell _2,\ell _3).$ 
The classical construction of a line spread comes from a pencil of hyperbolic quadrics
(see e.g. \cite{MR840877}, Lemma 17.1.1) and  
it has the property that if it contains any three lines of a regulus
${\cal R}(\ell _1, \ell _2,\ell _3),$  then it contains each 
of the $q+1$ lines of ${\cal R}(\ell _1, \ell _2,\ell _3).$ 
This type of spread is called \emph{regular}.
A line spread in $\mathrm{PG}(3,q)$ is called \emph{aregular}, if it 
contains no regulus. An example of an aregular spread can be found 
in \cite{MR840877}, Lemma 17.3.3.

\section{Factorizations arising from affine spaces}
\label{Affine}

In this section, we investigate the spreads and packings 
of $\mathrm{AG}(n,q)$ and the corresponding factorizations of
multigraphs. In each case we apply Lemma \ref{basic}, so we identify
the points of $\mathrm{AG}(n,q)$ with the vertices of the complete
multigraph.    

\begin{theorem}
\label{affgen}
Let $q$ be a prime power, $i<n$ be positive integers and 
$\lambda_i=\genfrac[]{0pt}{2}{n-1}{i-1}_q.$
Then there exists a simple ($q^i-1$)-factorization ${\cal F}^i$ of 
$\lambda_i K_{q^{n}}.$
${\cal F}^i$ is decomposable if and only if there exists an $f$-fold 
$(i-1)$-spread in
$\mathrm{PG}^{(i-1)}(n-1,q)$ for some $1\leq f<\lambda_i .$
\end{theorem}

{\bf Proof.}  
Consider the $n$-dimensional affine space as 
$\mathrm{AG}(n,q)=\mathrm{PG}(n,q)\setminus {\mathcal H}_{\infty }$
where ${\mathcal H}_{\infty }$ is isomorphic to $\mathrm{PG}(n-1,q).$ 
Take the design $D=\mathrm{AG}^{(i)}(n,q)$ and apply Lemma \ref{basic}.
If $\Pi ^{i-1}_j$ is an $(i-1)$-space of ${\mathcal H}_{\infty },$ then
the set of the $q^{n-i}$ parallel affine $i$-spaces through $\Pi ^{i-1}_j$ 
is an $i$-spread of $D.$ 
This spread induces a $(q^i-1)$-factor $F^i_j$ for $j\in\{1,\dots,r\}.$ 
If $\Pi ^{i-1}_1,\Pi ^{i-1}_2,\dots ,\Pi ^{i-1}_g$ are
distinct $(i-1)$-spaces of ${\mathcal H}_{\infty }$ and they form 
an $f$-fold spread, then $f=(g(q^i-1))/(q^n-1),$ and the union of the 
corresponding $(q^i-1)$-factors $F^i_j,$ for $j=1,2,\dots ,g,$ 
gives a ($q^i-1$)-factorization of $fK_{q^{n}}.$
Distinct $(i-1)$-spaces of ${\mathcal H}_{\infty }$ obviously
define distinct $(q^i-1)$-factors, so this factorization is simple.
In particular if we consider all $(i-1)$-spaces of 
${\mathcal H}_{\infty },$
then 
$$g=\genfrac[]{0pt}{0}{n}{i}_q, \, \, f=\genfrac[]{0pt}{0}{n}{i}_q
\frac{q^i-1}{q^n-1}=\genfrac[]{0pt}{0}{n-1}{i-1}_q=\lambda _i,$$
hence the union of the corresponding factors
gives a simple ($q^i-1$)-factorization ${\cal F}^i$ of $\lambda_i K_{q^{n}}.$ 

Suppose that ${\cal F}^i$ is decomposable, 
then there exist two positive integers
$\mu _1$ and $\mu _2$ such that
$\mu _1+\mu _2=\lambda _i$ and ${\cal F}^i$ can be written
as the union  ${\cal F}^i= {\cal F}_1\cup {\cal F}_2$;  ${\cal F}_1$ and
${\cal F}_2$ are $(q^i-1)$-factorizations of $\mu _1K_{q^{n}}$ and
$\mu _2K_{q^{n}},$ respectively, having no $(q^i-1)$-factors in common,
since ${\cal F}^i$ is simple.
For $h=1,2,$ the relation 
$\mu _h{{q^n}\choose 2}={{q^i}\choose 2}q^{n-i}|{\cal F}_h|$ holds, whence
$\mu _h(q^n-1)=(q^i-1)|{\cal F}_h|.$
Without loss of generality we can set 
${\cal F}_1= \cup _{j=1}^{f_1}F^i_j$
with $f_1=(\mu _1(q^n-1))/(q^i-1),$ and   
${\cal F}_2= {\cal F}^i\setminus {\cal F}_1,$
$f_2=|{\cal F}_2|.$

Let $u_1$ and $u_2$ be two affine points and let $w$ be 
the point at infinity of the line $u_1u_2$. 
Since ${\cal F}_h$ is a factorization of 
$\mu _hK_{q^{n}},$ there are exactly $\mu _h$
factors of ${\cal F}_h$ containing the edge $[u_1,u_2],$ say
$F^i_{j_1},F^i_{j_2},\dots ,F^i_{j_{\mu _h}}.$
The edge $[u_1,u_2]$ belongs 
the $F^i_{j_s}$ if and only if $w\in \Pi ^{i-1}_{j_s}$ for 
every $1\leq s\leq \mu _h.$ 
This happens if and only if $\cup _{j=1}^{f_h}\Pi ^{i-1}_j$ contains 
each point of ${\mathcal H}_{\infty }$  exactly $\mu _h$ times, which means 
that $\cup _{j=1}^{f_h}\Pi ^{i-1}_j$ is a $\mu _h$-fold spread in 
${\mathcal H}_{\infty },$ for every $h=1,2.$ It is thus proved
that if ${\cal F}^i$ is decomposable, then $\mathrm{PG}^{(i-1)}(n-1,q)$
posesses an $f$-fold spread for some $1\leq f< \lambda _i.$

Vice versa, suppose that there exists a $\mu _1$-fold spread in 
$\mathrm{PG}^{(i-1)}(n-1,q)$ for some $1\leq \mu _1< \lambda _i.$
Let ${\cal F}_1= \cup _{j=1}^{f_1}F^i_j$ be a $\mu _1$-fold 
spread in ${\mathcal H}_{\infty }.$ Then $|{\mathcal F}_1|=f_1=
\mu _1(q^n-1)/(q^i-1).$ Let ${\mathcal T}$  be the set of all 
$(i-1)$-dimensional subspaces in ${\mathcal H}_{\infty }$ and let
${\mathcal F}_2={\mathcal T}\setminus {\mathcal F}_1.$ Then 
$|{\mathcal T}|=\genfrac[]{0pt}{1}{n}{i}_q,$ hence  
$$|{\mathcal F}_2|=\genfrac[]{0pt}{0}{n}{i}_q-\mu _1(q^n-1)/(q^i-1)=
\left( \genfrac[]{0pt}{0}{n-1}{i-1}_q-\mu _1\right) \frac{q^n-1}{q^i-1},$$
so if 
$\mu _2= \genfrac[]{0pt}{1}{n-1}{i-1}_q-\mu _1,$ then
${\mathcal F}_2$ is a $\mu _2$-fold spread in ${\mathcal H}_{\infty }$
and $1\leq \mu _2<\lambda _i$ holds.   

As we have already seen, ${\mathcal F}_h$ defines  
a ($q^i-1$)-factorization of $\mu _hK_{q^{n}}$ for $h=1,2.$
Then ${\cal F}^i= {\cal F}_1\cup {\cal F}_2,$ 
because $\mu _1+\mu _2=\lambda _i.$ Hence the 
($q^i-1$)-factorization  ${\cal F}^i$ of $\lambda _iK_{q^{n}}$
is decomposable.
\Qed 

\begin{corollary}
If $\gcd (i,n)>1$ then the $(q^i-1)$-factorization ${\cal F}^i$ 
of $\lambda _iK_{q^{n}}$ is decomposable.
\end{corollary}

\noindent
{\bf Proof.}
Let $1<l+1$ be a divisor of $\gcd (i,n).$ Then 
it follows from Theorem \ref{spreadexists} that there exists
an $(q^{i}-1)/(q^{l+1}-1)$-fold spread in ${\mathcal H}_{\infty },$
so ${\cal F}^i$ is decomposable. 
\Qed

To decide the decomposability of ${\cal F}^i$ in the cases $\gcd (i,n)=1$ 
is a hard problem in general. We prove its indecomposability in the following
important case.

\begin{theorem}
The $(q^{n-1}-1)$-factorization ${\cal F}^{n-1}$ of 
$(q^{n-1}-1)/(q-1)K_{q^{n}}$ is indecomposable.
\end{theorem}

\noindent
{\bf Proof.}
It is enough to prove that if  
$\cup _{j=1}^{g}\Pi ^{n-2}_j$ is an $f$-fold $(n-2)$-spread in 
${\mathcal H}_{\infty },$ then $\cup _{j=1}^{g}\Pi ^{n-2}_j$
consists of all $(n-2)$-dimensional subspaces of ${\mathcal H}_{\infty },$ 
because this implies $f=\lambda _{n-1},$ so the statement follows from
Theorem \ref{affgen}.
 
Each $\Pi ^{n-2}_j$ contains exactly $(q^{n-1}-1)/(q-1)$ points, thus the 
standard double counting of the point-subspace pairs $p\in \Pi ^{n-2}_j$ in 
${\mathcal H}_{\infty }$ gives 
\[g \frac{q^{n-1}-1}{q-1}=f\frac{q^n-1}{q-1},\] hence 
\[f=\frac{g (q^{n-1}-1)}{q^n-1}.\] 
But $\mathrm{gcd}(q^n-1,q^{n-1}-1)=q-1$ and $f$ is an integer,
so $g \geq (q^n-1)/(q-1)$ which implies
$g =(q^n-1)/(q-1),$ hence $f=\lambda _{n-1}.$ 
\Qed 

\noindent
In particular if $n=2,$ we get the following.
 
\begin{corollary} 
\label{coro3.4}
If $q$ is a prime power then there exists a simple and indecomposable
$(q-1)$-factorization of $K_{q^2}$.
\end{corollary}

If $q=2^r$ then each $(q^i-1)$-factor in ${\cal F}^{i}$ is the vertex-disjoint
union of $2^{r-i}$ complete graphs on $2^i$ vertices.
It is well-known that these graphs can be partitioned 
into one-factors in many ways (but not in all the ways, it was proved by
Hartman and Rosa \cite{hr}, that there is no cyclic one-factorization
of $K_{2^i}$ for $i\geq 3$), 
hence Theorem \ref{affgen} implies several one-factorizations of 
$\lambda _iK_{2^r}.$ 

Each of the one-factorizations arising from 
${\cal F}^i$ is simple, because distinct $(i-1)$-dimensional subspaces   
define distinct $(q^i-1)$-factors of ${\cal F}^i,$ and the one-factors 
of $\lambda _iK_{q^n}$ arising from distinct $(q^i-1)$-factors of ${\cal F}^i$
are distinct, because they are the union of $q^{n-i}$ one-factors on 
$q^i$ vertices of a connected component.

There are both decomposable and indecomposable  
one-factorizations among these examples.
We show it in the smallest case $q=2, \,n=3.$ Let ${\cal F}^2$ be the
$3$-factorization of $3K_8$ induced by $\mathrm{AG}(3,2).$ 

Let $\mathrm{PG}(3,2)=\mathrm{AG}(3,2)\cup \mathcal{H}_{\infty }.$ 
Then $\mathcal{H}_{\infty}$ is
isomorphic to the Fano plane. Let its points be $0,1,2,3,4,5$ and $6$ such
that for $j=0,1,\dots ,6,$ the triples $L_j=(j,j+1,j+3)$  
form the lines of the plane, where the addition is taken modulo 7. Now the
$3$-factors of ${\cal F}^2$ can be described in the following way. Let $a$ be a
fixed point in $\mathrm{AG}(3,2).$ Then $L_j$ defines a $3$-factor $F^2_j$
whose connected components are complete graphs $K_{2^i}=K_4.$ 
Let $L_{j,a}$ be the complete graph containing $a,$ and let $L_{j,\overline{a}}$
be the other component of $F^2_j.$

$\mathcal{H}_{\infty}$ defines one-factors and a one-factorization 
of $K_8$ in the
following obvious way. The edge joining two points of $\mathrm{AG}(3,2),$ say
$b$ and $c,$ belong to the one-factor $G_s$ if and only if $b,c$ and $s$ are
collinear points in $\mathrm{PG}(3,2).$ Then 
${\cal G}=\cup _{s=0}^6G_s$ is a one-factorization of $K_8.$

We can define a decomposable one-factorization of $3K_8$ in the 
following way. Take
$L_{j,a}$ and $L_{j,\overline{a}}$ and let $s\in L_j$ be any point. Then $G_s$
gives a one-factor of $L_{j,a}$ and a one-factor of $L_{j,\overline{a}}.$ 
Hence ${\cal G}_j=\cup _{s\in L_j}G_s$ 
is the union of three one-factors of $3K_8,$ and
${\cal G}'=\cup _{j=0}^6{\cal G}_j$ is a one-factorization of $3K_8.$ 

In $\mathcal{H}_{\infty}$ there are three lines
through the point $s,$ hence ${\cal G}'$  
contains each one-factor $G_s$ three times.
Thus ${\cal G}'$ is decomposable, because it is obviously the
union of three copies of ${\cal G}.$

But we can define an indecomposable one-factorization, too. Let
$L_j$ be a line in $\mathcal{H}_{\infty},$
take $L_{j,a}$ and $L_{j,\overline{a}}$ and let
$M_{j}^1$ be the one-factor which contains the following pairs of points 
in $\mathrm{AG}(3,2):$ 

-- $(b,c)$ if $b,c\in L_{j,a}$ and $b,c,j$ are collinear 
in $\mathrm{PG}(3,2).$ 

-- $(b,c)$ if $b,c\in L_{j,\overline{a}}$ and $b,c,j+1$ are collinear 
in $\mathrm{PG}(3,2).$ 

Let $M_{j}^2$ be the one-factor which contains the following pairs of points in $\mathrm{AG}(3,2):$ 

-- $(b,c)$ if $b,c\in L_{j,a}$ and $b,c,j+1$ are collinear 
in $\mathrm{PG}(3,2).$ 

-- $(b,c)$ if $b,c\in L_{j,\overline{a}}$ and $b,c,j+3$ are collinear 
in $\mathrm{PG}(3,2).$

Finally let $M_{j}^3$ be the one-factor which contains the 
following pairs of points in $\mathrm{AG}(3,2):$ 

-- $(b,c)$ if $b,c\in L_{j,a}$ and $b,c,j+3$ are collinear in $\mathrm{PG}(3,2).$ 

-- $(b,c)$ if $b,c\in L_{j,\overline{a}}$ and $b,c,j$ are collinear in $\mathrm{PG}(3,2).$ 

Then ${\cal M}_{j}=\cup _{t=1}^3M_{j}^t$ is a union of three 
one-factors of $3K_8,$ and ${\cal M}=\cup _{j=0}^6{\cal M}_{j}$ 
is a one-factorization of $3K_8.$ 

Suppose that this one-factorization is decomposable. Then it contains a one-factorization 
${\cal E}$ of $K_8.$ ${\cal E}$ is the union of seven one-factors. 
We may assume without loss of generality, that $M_{0}^1$ belongs to ${\cal E}$. 
It contains an edge through $a,$ let it be $(a,b),$ and a pair $(c,d)$ for which the lines 
$ab$ and $cd$ are parallel lines in $\mathrm{AG}(3,2).$ 
There are two more lines in the parallel class of $ab,$ say $ef$ and $gh.$ 
It follows from the definition of the one-factors that exactly one of them contains 
the pairs $(e,f)$ and $(a,b),$ another one contains the pairs  $(e,f)$ and 
$(c,d),$ and a third one contains the pairs $(e,f)$ and $(g,h).$ 
But ${\cal E}$ contains each pair exactly once, hence it must contain the
one-factor containing the pairs $(e,f)$ and $(g,h).$ But this is a one-factor of
type $M_{0}^t,$ where $t\neq 1.$ Hence ${\cal E}$ contains $M_{0}^t$ 
where $t=2$ or 3. If we repeat the previous argument, we get that ${\cal E}$ 
must contain  $M_{0}^l$ for $1\neq l\neq t,$ too. 
Thus ${\cal E}$ is the union of triples of type $M_{j}^t,$ $t=1,2,3,$ but
this is a contradiction, because ${\cal E}$ consists of seven 
one-factors.

\section{Factorizations arising from projective spaces}
\label{Proj}

There are two basic types of partitioning the point-set of 
finite projective spaces. Both types give factorizations of
some multigraphs. In this section we discuss these constructions.

\subsection{Spreads consisting of subspaces} 
\label{1-spreads}
It is easy to construct 
spreads in $\mathrm{PG}^{(i)}(n,q),$ Theorem \ref{spreadexists}
gives a necessary and sufficient existence condition. Packings are 
much more complicated objects. Only a few packings in $\mathrm{PG}^{(1)}(n,q)$
have been constructed so far. In each case of 
the known packings either $n$ or $q$ 
satisfies some conditions.

\begin{theorem}[Beutelspacher, \cite{MR0341270}] 
\label{Beutelspacher}
Let $1<k$ be an integer and let $n=2^{k}-1$.
Then there exists a packing in $\mathrm{PG}^{(1)}(n,q).$  
\end{theorem}

\begin{theorem}[Baker, \cite{MR0416937}] 
\label{Baker}
Let $1<k$ be an integer. Then
there exists a packing in $\mathrm{PG}^{(1)}(2k-1,2)$.
\end{theorem}

\noindent
Applying the Basic Construction Lemma, we get the following 
existence theorems. 

\begin{corollary}
\label{c1}
Let $q$ be a prime power, $1<k$ be an integer and 
$v=\tfrac{q^{2^k}-1}{q-1}.$ Then there exists a $q$-factorization of $K_v$ 
induced by a line-packing in $\mathrm{PG}(2^k-1,q)$.
\end{corollary}

\begin{corollary}
\label{c2}
Let $1<k$ be an integer and $v=\tfrac{q^{2k}-1}{q-1}.$
There exists a $2$-factorization $K_v$  
induced by a line-packing in $\mathrm{PG}(2k-1,2)$.
\end{corollary}

If $k=2$ then Corollary \ref{c2} gives a solution 
of {\it Kirkman's fifteen schoolgirls problem}, which was first posed in 1850
(for the history of the problem we refer to \cite{biggs}), while 
Corollary \ref{c1} gives a solution of the generalised problem in the case
of $(q^2+1)(q+1)$ schoolgirls. 

The complete classification of packings in $\mathrm{PG}^{(i)}(n,q)$ is known 
only in the case $i=1,\, n=3$ and $q=2.$ There are 240 projectively distinct
packings of lines in $\mathrm{PG}(3,2)$ (see \cite{MR840877}, Subsection 17.5).

If $\gcd (q+1,3)=3$, then
there is a construction of spreads in $\mathrm{PG}^{(1)}(3,q)$ due to 
Bruen and Hirschfeld \cite{brh} which is completly different from the 
constructions of Theorems \ref{Beutelspacher} and \ref{Baker}.
It is based on the geometric properties of twisted cubics.

A normal rational curve of order 3 in $\mathrm{PG}(3,q)$ 
is called \emph{twisted cubic}. 
It is known that a twisted cubic is projectively equivalent to the set of
points $\{ (t^3:t^2:t:1)\, \colon t\in \mathrm{GF}(q)\} \cup \{ (1:0:0:0)\} .$
In \cite{MR840877} it was shown that there exist aregular spreads given
by a twisted cubic. For a detailed description of twisted cubics and the proofs
of the following theorems we refer to \cite{MR840877}, Section 21.

\begin{theorem}
\label{Hirsschfeld1}
Let $G_q$ be the group of projectivities in $\mathrm{PG}(3,q)$ fixing a
twisted cubic $\mathcal{C} .$ Then 
\begin{itemize}
\item
$G_q \cong \mathrm{PGL}(2,q) $ and it acts
triply transitively on the points of $\mathcal{C}.$ 
\item
If $q \geq 5$ then the number of twisted cubics in $\mathrm{PG}(3,q)$ is
$q^5(q^4-1)(q^3-1).$
\end{itemize} 
\end{theorem}

\begin{theorem}
Let $\mathcal{C}$ be a twisted cubic in $\mathrm{PG}(3,q).$
If $\gcd(q+1,3)=3$, then there exists a spread 
in $\mathrm{PG}^{(1)}(3,q)$ induced by $\mathcal{C}$.
\end{theorem}

\noindent
Using the spreads associated to twisted cubics and the Basic Construction 
Lemma, we get the following multigraph factorization.

\begin{theorem}
Let  $q\geq 5$ be a prime power, 
$\lambda =q^5(q^4-1)(q-1)$ and $v={q^3+q^2+q+1}.$
If $\gcd (q+1,3)=3,$ then there exists a simple
$q$-factorization of $\lambda K_v$ induced by the set of 
twisted cubics in $\mathrm{PG}(3,q)$.
\end{theorem}

\noindent
{\bf Proof.} 
Let $\mathsf{C} $ be the set of twisted cubics in $\mathrm{PG}(3,q)$.
For $\mathcal{C}\in \mathsf{C}$ let $\mathcal{L}_{\mathcal{C}} $ be the spread 
in $\mathrm{PG}^{(1)}(3,q)$ induced by $\mathcal{C}$.
If $\ell $ is a line and $c_{\ell }$ denotes the number of
twisted cubics $\mathcal{C}$ with the property that $\ell $ belongs to $\mathcal{L}_{\mathcal{C}},$
then it follows from Theorem \ref{Hirsschfeld1} that $c_{\ell }$ does not
depend on $\ell .$ Hence
$$
\begin{aligned}
c_{\ell} & = \frac{|\{ \textrm{twisted cubics in }\mathrm{PG}(3,q)\} |
\times|\{ \textrm{lines in a spread of }\mathrm{PG}(3,q)\} |}
{|\{ \textrm{lines in }\mathrm{PG}(3,q)\} | } \\
 & =\frac{q^{5}(q^{4}-1)(q^{3}-1)\times(q^{2}+1)}{(q^{2}+1)(q^{2}+q+1)}=q^{5}(q^{4}-1)(q-1).
\end{aligned}
$$

Thus $\mathsf{C} $ induces a $|\mathsf{C}|$-fold spread in 
$\mathrm{PG}^{(1)}(3,q)$. 
Each spread $\mathcal{L}_{\mathcal{C}}$ induces a $q$-factor in $K_v$, 
therefore the Basic Construction Lemma gives that
$\underset{\mathcal{C}\in\mathsf{C}}{\bigcup}\mathcal{L}_{\mathcal{C}}$ 
is a $q$-factorization of $\lambda K_v$.
Any two distinct twisted cubics define different spreads, hence the
factorization is simple by definition.
\Qed

\subsection{Constructions from subgeometries}

If the order of the base field is not prime, then projective spaces 
can be partitioned by subgeometries. Let $1<k$ be an integer.
Since $\mathrm{GF}(q)$ is a subfield of $\mathrm{GF}(q^k),$ so
$\mathrm{PG}(n,q)$ is naturally embedded into $\mathrm{PG}(n,q^k)$ if
the coordinate system is fixed. Any $\mathrm{PG}(n,q)$ embedded into 
$\mathrm{PG}(n,q^k)$ is called a \emph{subgeometry}. 
Using cyclic projectivities one can prove that any $\mathrm{PG}(n,q^k)$ 
can be partitioned by subgeometries $\mathrm{PG}(n,q)$. 
For a detailed description of cyclic projectivities, subgeometries, and the 
proofs of the following three theorems we refer to \cite{MR554919}, Section 4.

\begin{theorem}[\cite{MR554919}, Lemma 4.20]
\label{theo21}
Let $s(n,q,q^k)$ denote the number of subgeometries $\mathrm{PG}(n,q)$ in 
$\mathrm{PG}(n,q^k).$ Then 
$$s(n,q,q^k)=q^{\binom{n+1}{2}(k-1)}\overset{n+1}{\underset{i=2}{\prod}}\tfrac{q^{ki}-1}{q^{i}-1}.$$
\end{theorem}

\begin{theorem}[\cite{MR554919}, Theorem 4.29]
$\mathrm{PG}(n,q^k)$ can be partitioned into 
$\theta (n,q,q^k) = \tfrac{(q^{k(n+1)}-1)(q-1)}{(q^{k}-1)(q^{n+1}-1)}$ 
disjoint subgeometries $\mathrm{PG}(n,q)$ if and only if $\gcd(k,n+1)=1$.
\end{theorem}

\begin{theorem}[\cite{MR554919}, Theorem 4.35]
Suppose that $\gcd(k,n+1)=1.$ 
Let $p_0 (n,q,q^k)$ denote the number of projectivities
which act cyclically on a $\mathrm{PG}(n,q)$ of
$\mathrm{PG}(n,q^k)$ such that determine different partitions. Then 
$$p_0 (n,q,q^k)=
q^{k \binom{n+1}{2}}\frac{\overset{n}
{\underset{i=1}{\prod}}(q^{ki}-1)}{n+1}.$$ 
Any given subgeometry $\mathrm{PG}(n,q)$ is contained in 
$$\rho_0(n,q)
=q^{\binom{n+1}{2}}\frac{\overset{n}{\underset{i=1}{\prod}}(q^i-1)}{n+1}$$
of these partitions.
\end{theorem}

We can consider the partitions of the point-set of 
$\mathrm{PG}(n,q^k)$ by subgeometries $\mathrm{PG}(n,q)$.

Each partition of $\mathrm{PG}(n,q^k)$
into subgeometries $\mathrm{PG}(n,q)$
defines a $\left(\frac{q(q^n-1)}{q-1}\right) $-factor of $K_v,$
with $v=\frac{q^{k(n+1)}-1}{q^k-1}.$ Each projectivity which acts 
cyclically on a $\mathrm{PG}(n,q)$ defines a   
$\left(\frac{q(q^n-1)}{q-1}\right) $-factorizations of the 
corresponding complete multigraph.

\begin{theorem}
\label{subgeoconstruct}
Let $q$ be a prime power, $1<k$ and $n$ be positive integers for which 
$\gcd(k,n+1)=1$ holds. Let 
$\lambda =\tfrac{q^{\binom{n+1}{2}k}(q^{k}-1)(q^n-1)}{q^{k-1}(n+1)(q-1)}
\overset{n-1}{\underset{i=1}{\prod}}(q^{ki}-1)$ 
and $v=\tfrac{q^{k(n+1)}-1}{q^{k}-1}.$
Then there exist a simple $\left( q\frac{q^n-1}{q-1}\right) $-factorization of 
$ \lambda K_v$ 
induced by the set of those projectivities which act cyclically on a
$\mathrm{PG}(n,q)$ of $\mathrm{PG}(n,q^k)$ such that they determine different 
partitions.
\end{theorem}

\noindent
{\bf Proof.} 
It follows from Theorem \ref{theo21} that the number $S_e$ of subgeometries 
$\mathrm{PG}(n,q)$ through two points of $\mathrm{PG}(n,q^k)$ is
$$
\begin{aligned}
S_e & = \frac{s(n,q,q^k)\times\left|\{ \textrm{points in }\mathrm{PG}(n,q) \} 
\right|\times(\left|\{ \textrm{points in }\mathrm{PG}(n,q)\} \right|-1)}
{\left|\{ \textrm{points in }\mathrm{PG}(n,q^k)\} \right|\times
(\left|\{\textrm{points in }\mathrm{PG}(n,q^k)\} \right|-1)} \\
 & =\frac{q^{\binom{n+1}{2}(k-1)}(q^{k}-1)}{q^{k-1}(q-1)}\overset{n-1}{\underset{i=1}{\prod}}\frac{q^{ki}-1}{q^{i}-1}.
\end{aligned} 
$$
Each cyclic projectivity determines different partitions, hence it 
determines different factors. 
Thus $\lambda =S_e \times \rho_0 (n,q) $.
\Qed

We cannot decide the decomposability of the factorization construted in 
the previous theorem in general, but we can prove the existence of 
indecomposable factorizations in 
some cases. To do this we need the following result from number theory.

\begin{lemma}[\cite{MR554919}, Lemma 4.24]
\label{divides}
If $r$, $s$ and $x$ are positive integers with $x>1$, then $\tfrac{(x^{rs}-1)(x-1)}{(x^{r}-1)(x^{s}-1)}$ 
is an integer if and only if $\gcd(r,s)=1$.
\end{lemma}

\noindent
We apply it in a particular case.

\begin{proposition}
Let $q$ be a prime power, $1<k$ and $n$ be positive integers for which 
$\gcd(k,n+1)=1$ and $\gcd (k,n)\neq 1$ hold. Let  
$d=\gcd \left( \tfrac{q^{kn}-1}{q^{k}-1},\tfrac{q^{n}-1}{q-1}\right) ,$
$v=\tfrac{q^{k(n+1)}-1}{q^{k}-1}$ and $m=q\frac{q^n-1}{q-1}.$ 
Suppose that $\cal F$ is   
an $m$-factorization of $ \lambda K_v$  
for some $\lambda $ such that each
factor is the disjoint union of 
$\theta (n,q,q^k)$ complete graphs on $(q^{n+1}-1)/(q-1)$ vertices. 
If $f$ denotes the number of $m$-factors in $\cal F$ 
then $\tfrac{q^{n}-1}{d(q-1)}$
divides $\lambda$ and $q^{k-1}\tfrac{q^{kn}-1}{d(q^{k}-1)}$ divides $f.$
\end{proposition}

\noindent
{\bf Proof.} 
The standard double counting gives 
\[\lambda\times\tbinom{v}{2}=\tbinom{m+1}{2}\times\theta(n,q,q^{k})\times f,\]
thus $\lambda\times
q^{k-1}\tfrac{q^{kn}-1}{d(q^{k}-1)}=f\times\tfrac{q^{n}-1}{d(q-1)}$. Because 
of Lemma \ref{divides}, $\tfrac{q^{n}-1}{d(q-1)}$ divides $\lambda ,$ hence
$q^{k-1}\tfrac{q^{kn}-1}{d(q^{k}-1)}$ 
divides $f$.
\Qed 

As a direct corollary of the previous proposition we get the following result 
about the indecomposibility of the factorizations constructed 
in Theorem \ref{subgeoconstruct} .
\begin{theorem}
Let $q$ be a prime power, $1<k$ and $n$ be positive integers for which 
$\gcd(k,n+1)=1$ and $\gcd (k,n)\neq 1$ hold. Let  
$d=\gcd \left( \tfrac{q^{kn}-1}{q^{k}-1},\tfrac{q^{n}-1}{q-1}\right) ,$
$v=\tfrac{q^{k(n+1)}-1}{q^{k}-1}$ and $m=q\frac{q^n-1}{q-1}.$
Then there exist a simple and indecomposable
$m$-factorization of $ \lambda K_v,$ where  
$\lambda = t\tfrac{q^{n}-1}{d(q-1)}$ for
some $t$ in
$\{1,\dots,d\tfrac{q^{\binom{n+1}{2}k}(q^{k}-1)}{q^{k-1}(n+1)}\overset{n-1}{\underset{i=1}{\prod}}(q^{ki}-1)\}.$
\end{theorem}

\bigskip

\bigskip
\noindent
{\bf \Large{Acknowledgement}}

\smallskip
\noindent
The authors are grateful to the anonymous reviewers for their
detailed and helpful comments and suggestions.

\begin{flushleft}
Gy\"{o}rgy Kiss \\
Department of Geometry and MTA-ELTE GAC Research Group \\
E\"{o}tv\"{o}s Lor\'{a}nd University \\
1117 Budapest, P\'{a}zm\'{a}ny s. 1/c, Hungary \\
e-mail: {\sf kissgy@cs.elte.hu}
\end{flushleft}

\begin{flushleft}
Christian Rubio-Montiel \\
Instituto de Matem{\'a}ticas \\ 
Universidad Nacional Aut{\'o}noma de M{\'e}xico \\ 
Ciudad Universitaria, 04510, D.F., Mexico \\
e-mail: {\sf christian@matem.unam.mx}
\end{flushleft}

\end{document}